\input amstex
\define\aut{\operatorname{Aut}}
\define\auto{\operatorname{Aut^0}}
\define\supp{\operatorname{supp}}
\define\csupp{\operatorname{csupp}}
\define\lk{\operatorname{lk}}
\define\rk{\operatorname{rk}}
\define\red{\overset\text{red}\to=}
\define\ett{\operatorname{and}}
\documentstyle{amsppt}
\rightheadtext{A Generating Set for $\aut G$}
\leftheadtext{L.J. Corredor $\ett$ M.A. Gutierrez}
\TagsOnRight
\topmatter
\title
A Generating Set for the Automorphism Group of \linebreak a Graph Product of Abelian Groups
\endtitle
\author
L.J. Corredor and M.A. Gutierrez
\endauthor
\address
Departamento de Matem\'aticas,  Universidad de los Andes,  Bogot\'a, Colombia
\linebreak
and Mathematics Department,  Tufts University,  Medford, MA 02155
\endaddress

\abstract
We find a set of generators for the automorphism group $\aut G$ of a graph product $G$ of finitely generated abelian groups entirely from a certain labeled graph. In addition, we find generators for the important subgroup $\aut^\star G$. We follow closely the plan of M. Laurence's paper [11].
\endabstract
\email
LJC: lcorredo\@uniandes.edu.co
MAG: mauricio.gutierrez\@tufts.edu
\endemail
\endtopmatter
\document
\head
1. Introduction
\endhead
A combinatorial graph $\Gamma$ is a pair $(V, E)$, where $V$ is a finite set of vertices and where $E$ is a set of edges which may be considered as subsets of $V$ of cardinality precisely equal to 2. Thus edges are unoriented, they have two distinct endpoints and two distinct vertices determine at most one edge.

Let $\{G_v\}$ be a family of non-trivial groups indexed by $V$. Define the graph product $W(\Gamma, \{G_v\})$ to be the free product $\prod^*_{v\in V} G_v$ modulo the subgroup generated by the $[G_x , G_y]$ for $\{x,y\}\in E$.

One way of defining a family $\{G_v\}$ is by letting $P\subset \Bbb N\cup\{\infty\}$ be the set of all non-trivial prime powers and $\infty$ and defining a labeled graph $(\Gamma, o)$ as a pair, where $o:V\to P$. Then we take $G_v$ to be the cyclic group of order $o(v)$ generated by $v\in V$. Thus $G_v$ is infinite cyclic if $o(v)=\infty$ and finite cyclic of order $o(v)$ otherwise. Then we denote the graph product by $W(\Gamma, o)$ or simply by $W(\Gamma)$.

A graph isomorphism $\gamma:\Gamma\to\Gamma'$ is a labeled graph isomorphism $(\Gamma,o)\to(\Gamma',o')$ if  $o'(\gamma(v))=o(v)$.

Suppose $\Delta=(X,F)$ is a graph and $\{H_x\}$ is a family of finitely generated abelian groups indexed by $X$. In [6] it is shown that if $G\simeq W(\Delta,\{H_x\})$, then
\roster
\item  there exists a labeled graph $(\Gamma,o)$ with $G\simeq W(\Gamma,o)$.
\item  such labeled graph is unique up to isomorphism.
\endroster

Thus the group $G$ determines $(\Gamma,o)$. Our purpose is to find a set of generators for $\aut G$ using the properties of $(\Gamma,o)$.

To define elements of $\aut G$, we usually profit from the fact that the vertices of $\Gamma$ generate $G$ and so we start with a map $\phi:V\to G$ and say that $\phi$ is well-defined (see [2], \S 2.3) if it extends to an endomorphism of $G$, which we also denote by $\phi$. 

Our list of generators comprises four families:
\roster
\item[1] Labeled graph automorphisms, already defined.
\item[2] Factor automorphisms (see [3], [4]) are defined by choosing a vertex $v$ and letting
$$\phi_v(z)=\cases
v^m &  z=v\\
z& z\neq v
\endcases
$$
we require that $m=-1$ if $o(v)=\infty$ and that $m\neq1$ and $(m,o(v))=1$ otherwise. As usual, $(m,o(v))=\gcd\{m,o(v)\}$.
\item[3] Dominated transvections. Suppose $x,y$ are two distinct vertices. Let
$$\tau_{x,y}(z)=
\cases
xy &\ \ \ z=x\\
z&\ \ \ z\neq x
\endcases
$$

In many cases (see Proposition 5.5) either $\tau_{x,y}$  or a power $\tau^q_{x,y}$ (which maps $x$ to $xy^q$) is well-defined. We choose a minimal $q\geq1$. 
\item[4] Partial conjugations. Again, choose a vertex $v$ and create a new graph $\Gamma\setminus v^\star$ by removing from $\Gamma$ the vertex $v$ and all its neighbors. This graph is not, in general, a connected graph. Let $K$ be a connected component of $\Gamma\setminus v^\star$. Define
$$\sigma_{K,v}(z)=
\cases
vzv^{-1}&\ \ \ z\in K\\
z&\ \ \ z\notin K
\endcases
$$
\endroster
The families 1, 2 and 4 are always well-defined. All are automorphisms when well-defined.

We now summarize our results: given $G$, a graph product of finitely generated abelian groups, it is possible  to write $G\simeq W(\Gamma, o)$. From $(\Gamma, o)$ construct the finite set $\Cal G$ consisting of all labeled graph automorphisms, factor automorphisms, well-defined dominated transvections and partial conjugations. Then we have
\proclaim{Main Theorem} The set $\Cal G$ generates $\aut G$.
\endproclaim
We devote a final small section to a subset $\Cal G^\star$ of $\Cal G$ that generates an important subgroup $\aut^\star G$ of $\aut G$ defined in [7]. In addition, we discuss the possibility of extending the results of [2] and finding a finite
presentation for $\aut G$.
The reader will see immediately that we  follow the plan in [10] and [11]; we quote many of the results without proof but on occasion we add our own proofs when the presence of
finite groups requires it, when they are unavailable elsewhere or when the original proofs are, at least to the authors, obscure. 

Conversations with X. Caicedo were very useful.

\head
2. Graphs
\endhead
If $\Gamma=(V,E)$ is a graph, a subgraph $\Sigma$ is a pair $(X,F)$ with $X\subset V$ and $F$ consists of the elements $\{x,y\}$ of $E$ with $x,y\in X$. Thus $X$ determines $\Sigma$ and we write $\Sigma= \langle X\rangle_{\Gamma}$. Often we mix subsets of $V$ with the subgraphs they generate and write, e.g., $v\in\Gamma$ to mean $v\in V$ and also use the expression $v$ for the subgraph whose only vertex is $v$. We write $|X|$ for the cardinality of the set $X$.

An important, and related graph is $\Gamma^c$, the complementary, or Coxeter, graph associated to $\Gamma$. It has the same vertices $V$ and $e$ is an edge of $\Gamma^c$ if and only if $e\notin E$. Here is a list of important subgraphs of $\Gamma$:
\roster
\item $\lk v=\langle\{x: \{v,x\}\in E\}\rangle_{\Gamma}$, the {\it link} of $v\in\Gamma$.
\item $v^\star=\langle\lk v \cup v\rangle_{\Gamma}$, the {\it star} of $v$.
\item $\Gamma\setminus v =\langle V\setminus v\rangle_{\Gamma}$ and $\Gamma\setminus v^\star=\langle V\setminus v^\star\rangle_{\Gamma}$.
\item if $\Sigma\subset\Gamma$, let $\lk\Sigma=\bigcap_{s\in\Sigma} \lk s$ and $\Sigma^\star=\bigcap_{s\in\Sigma}s^\star$.
\endroster.
\head
3. Graph Products
\endhead
We always put $\Gamma=(V,E)$ and $G=W(\Gamma,o)$. If $\Sigma\subset\Gamma$, $(\Sigma,o')$ is the labeled graph where $o'=o\restriction X$. By [7], Lemma 2.5, $W(\Sigma,o')$ embeds (as a retract) in $W(\Gamma,o)$. Then we write $W(\Sigma,o')=\langle\Sigma\rangle$ or even $\langle X\rangle$. The overlap in notation causes no ambiguity.   Let  $F(V)$ be the free group with basis $V$; then any element $g\in G$ can be written in at least one way as a word $a=v_1^{n_1}\cdots v_p^{n_p} $ of length $|a|:=\sum|n_i|$ and support $\supp a=\{v_1,\dots. v_p\}$. We say $a$ is reduced if it represents $g$ and it is of minimal length. Then we write $|g|:=|a|$ and $\supp g:=\supp a$. All this is well defined [6], \S2. We say that $g$ is cyclically reduced if, for any $v\in\Gamma$, $\supp v^{\pm1}g v^{\mp1}\supseteq\supp g$. In [14], Proposition II it is shown that each $g\in G$ has a unique reduced expression $g\overset\text{red}\to=wuw^{-1}$ with $u$ cyclically reduced (CR), and we refer to $u$ as the CR part of $g$. Then we write $\csupp g :=\supp u.$

If $g=v^{\pm1}g'$ reduced, we say that $v$ is a first letter of $g$. Last letters, initial segments and final segments are defined similarly. In [6] we prove
\proclaim{Rigidity Lemma}If $\theta\in {\aut G}$ then $\theta$ determines a graph automorphism $\gamma$ of $(\Gamma,o)$ which satisfies $\gamma (v)\in \csupp\theta(v)$ for all $v\in\Gamma$.
 \endproclaim

 \head
4. The Centralizer Theorem
\endhead
For an element $g\in G$ let $C(g)$ denote its centralizer subgroup. If $g\red wuw^{-1}$ then $C(g)=wC(u)w^{-1}$ and we may concentrate on the CR case only. Write $o(g)$ to indicate the order of $g$. Observe that  $o(g)$ does not necessarily lie in $P$. If $o(g)=\infty$, we write $g=r^n$, with $n\geq1$ maximal. If $o(g)<\infty$, we choose $n$ maximal in $1\leq n<o(g)$. Then we say $r=\sqrt g$ is a root of $g$.

Hereafter, $g$ is CR (so we use the symbol $u$ instead) and has support $S$ with $\Sigma=\langle S\rangle_{\Gamma}$. Then $\Sigma^c=\langle S\rangle_{\Gamma^c}$, so we write $\langle S\rangle^c$ for $\Sigma^c$.
\proclaim{Centralizer Theorem} (See [10], [14] and especially [1]) If $u$ and $\Sigma$ are as above
\roster
\item[1] If $\Sigma^c$ is connected, $u$ has a unique root.
\item[2] If $\Sigma^c=\Sigma_1^c\sqcup\cdots\sqcup\Sigma_n^c$ decomposes into $n$ distinct connected components, $\Sigma_i=\langle S_i\rangle^c$, then $u$ is a commuting product  $u=p_1\cdots p_n$, with $ \supp p_i=S_i$. Let $u_i=\sqrt p_i$; then $\supp u_i=S_i$ and we may write $u$ as a commuting product
$$
 u_1^{m_1}\cdots {u_n}^{m_n}.\tag4.1
$$

The centralizer of $u$ is

$$
C(u)= \langle{u_1}\rangle\times\cdots\langle{u_n}\rangle\times \langle\lk\Sigma\rangle.\tag4.2
$$
\endroster
\endproclaim
The product in (4.1) is called a {\it basic form} of $u$.

A general proof of this theorem is given in [1] but unfortunately its relation to $(\Gamma,o)$ is unclear and it does not immediately yield the formula in (4.2) that we will need in \S 5. Below we provide a sketch of an alternate, and elementary, proof using the ideas of [8] and [10]. Assume $v\in S$. If $\Sigma$ lies in some complete subgraph, the result is found in Lemma 2.2 of [6]. Assume then that $\Sigma$ does not lie in a complete subgraph. Then we can find $v\in\Sigma$ so that $\Sigma\nsubseteq v^\star$ and consequently, 
$$
  o(u)=\infty\ \ \ \text{and} \ \ |S|\geq2.\tag4.3
$$
\proclaim{Lemma 4.1} If $x$ is a last letter of $u$, $u\red u'x$ and $x^{-1}$ is not a first letter of $u'$, the equation $xu' \red y^{-1}u '' $, with $y$ not a last letter of $u''$,  implies that $x$ and $y$ commute, $x^{-1}$ is not a first letter of $u '' $ and $y$ is not a last letter of $u'$.
\endproclaim

This mostly follows from the fact that any two initial (final) letters of an element must commute.
\subhead
Reduction Types
 \endsubhead
 Assume that $w= x_n\cdots x_1$ ($x_i\in V^{\pm1}$) is reduced and that $u$ is CR. We analyze $wuw^{-1}$: define inductively
$$
u_0=u \quad u_i=x_i u^{i-1} x^{-1}_i.
$$

We shall write $u_{i-1}\overset{r}\to{\underset x_i\to\longrightarrow} u_i$ to mean that the reduction from $u_{i-1}\to u_i$ is of type $r$  and it involves the letter $x_i$. The Normal Form Theorem [8],
Theorem 3.9 and [10] Corollary 3.1.3, implies that there are five types $r$:
\roster
\item[0] $u_{i-1}\red x^{-1}_i u' x_i$. In particular, $u_{i-1}$ is not CR and $u_i=u'$.
\item[1] $x_i$  commutes with each vertex in $\supp u_{i-1}$. In particular $u_{i-1}=u_i$.
\item[2] $u_{i-1}\red x^{-1}_iu'$ and $x_i$ is not a last letter of $u'$. In particular, $u_i=u'x^{-1}_i$ is a cyclic permutation of $u_{i-1}$ with the same support. If $u_{i-1}$ is CR, so is $u_i$.
\item[3] $u_{i-1}\red u'x_i$ and $x_i^{-1}$ is not a first letter of $u'$. Same conclusions as in 2.
\item[4] $u_i\red x_i u_{i-1}x^{-1}_i$. In particular $u_i$ is not CR.
\endroster
\proclaim{Reduction Lemma} If $u,w$ are as above, then by rearranging the $x_i$ we may find numbers
$$
0=n_0\leq n_1\leq n_2\leq n_3\leq n_4=n \tag4.4
$$
such that, if $n_k<i\leq n_{k+1}$, $ k=0,1,2,3$, the reduction type $u_{i-1}\to u_i$ is $k+1$. Further, if $1\leq i\leq n_1<j\leq n_2<k\leq n_3$ then $x_i,\ x_j$ and $x_k$ commute in pairs.
\endproclaim
\demo{Proof} A complete proof is available in [10], Lemma 5.1 but unfortunately it contains an annoying (and correctable) gap. We run quickly over parts of the proof. If $u_0$ is CR then we show that there are no steps of type 0. as is  the case with the first step, and that the reductions of type 4 occur last. If $u_1,\dots,u_l$ involve only reduction types $\geq1$ and the next step is to be of type 0, then $u_l\red hu'h^{-1}$ and $u'$ is obtained from $u$ by reductions of types 1, 2 and 3 and so it is CR with $\supp u'=\supp u$. Say $h\red a_1\cdots a_t$; then $w\red w'h$. If $u_l\overset0\to{\underset x\to\longrightarrow}u_{l+1}$ then $u_l\red x^{-1}u'' x$ and $x\in\supp h$. If $x=a_k$, $x$ commutes with $a_1,\dots,a_{k-1}$ contradicting the fact that $h$ (and $w$) are reduced. Thus  it suffices to show that the initial reduction types occur in order. To simplify notation we use $\alpha, \beta$ etc instead if $u_{i-1}, u_i$ etc. If we have, say $\alpha\overset4\to{\underset x\to\longrightarrow}\beta\overset1\to{\underset y\to\longrightarrow}\gamma$ then $\beta=\gamma$ has support $\{x\}\cup\supp\alpha$. 
Since $y$ commutes with each element of $\supp\beta$, $y$ commutes with $x$ and with $\alpha$. Then we have $\alpha\overset1\to{\underset y\to\longrightarrow}\alpha\overset4\to{\underset x\to\longrightarrow}\gamma$ and a similar proof shows that you can put all reductions of type 1 first. Similarly reductions of type 4 commute with reductions of types 2 and 3 so the reductions of type 4 go last.

If  $\alpha$ is CR and $\alpha\overset3\to{\underset x\to\longrightarrow}\beta\overset2\to{\underset y\to\longrightarrow}\gamma$, then we are in the situation of Lemma 4.1: $\alpha\red u'x$ and $\beta=xu'=y^{-1}u''$, both CR so $x$ and $y$ commute and, 

$$
\alpha=y^{-1}x^{-1}u'' x\overset2\to{\underset y\to\longrightarrow}x^{-1}u''xy^{-1}=x^{-1}u'' y^{-1}x\overset3\to{\underset x\to\longrightarrow}u'' y^{-1}=\gamma,
$$
is the desired sequence. Only the expression for $\gamma$ is in reduced form.\qed 

 \enddemo
 \proclaim{Corollary 4.2}If $u$ is CR and $w\in G$, then $w$ has a reduced expression $w=w_4w_3w_2w_1$, with $|w_k|=n_k$ as in (4.4), where each letter of $w_k$ gives a reduction of type $k$, $w_1,\ w_2$ and $w_3$ commute in pairs and $wuw^{-1}$ is CR if and only if $w_4=1$.
 \endproclaim
 \proclaim{Proposition 4.3} If $u$ is CR, $w\in G$ commutes with $u$ and $\langle\supp u\cup\supp w\rangle^c$ is connected in $\Gamma^c$ then $w_1=w_4=1$, either $w_2$ or $w_3$ is trivial and $\supp w=\supp u$.
 \endproclaim
 \demo{Proof} Since $u=wuw^{-1}$ the latter is CR and so $w_4=1$ by Corollary 4.2. Then we prove sequentially
 \roster
\item $\supp w_2 \cup \supp w_3 \subset \supp u$,
 \item $\supp w\subset\supp u$,
 \item $\supp w=\supp u$,
 \item$w_1=1$ and
 \item either $w_2$ or $w_3$ is trivial.
 \endroster.

 Assertion (1) follows from $\supp u=\supp u_0=\supp u_k$, $1\leq k\leq n_3=n$. For (2), if $\langle\supp u\cup\supp w\rangle^c$ is connected and $\langle\supp w\rangle^c\setminus\langle\supp u\rangle^c\neq\emptyset$, then there are adjacent vertices $x\in\langle\supp w\rangle^c\setminus\langle\supp u\rangle^c$ and $y\in\langle\supp u\rangle^c$. Therefore $xy\neq yx$. But (1) implies $x\in\langle\supp w\rangle$ so $xy=yx$, a contradiction. For (3), write $w\red lw'l^{-1}$ with $w'$ CR which commutes with $l^{-1}ul$ and apply (2). Assertion (4) follows from the fact that (1), (2) and (3) imply $\supp u=\supp w_1\cup(\supp u\setminus\supp w_1)$. By connectedness one of the two terms is empty. If $\supp u=\supp w_1$, then it has more than 2 vertices, by equation (4.3), and the hypothesis of connectedness is contradicted. The last assertion is similar (cf. [10], Lemma 5.2 ).\qed
 \enddemo
 \proclaim{Lemma 4.4} Suppose $u$ is CR, $uw=wu$ and $\langle\supp u\rangle^c$ is connected. Then either $wu$ or $uw^{-1}$ is reduced as written.
 \endproclaim
See Lemma 5.3 in [10]. Compare to Proposition III.2 of [14].
 \subhead
 Free Products with Amalgamation
 \endsubhead
 In her thesis Green [8], Lemma 3.20 shows for a labeled graph $(\Gamma,o)$ and $v\in\Gamma$ that
 $$
 W(\Gamma)=W(v^\star)\ast_{W(\lk v)} W(\Gamma\setminus v).\tag4.5
 $$
 The hypothesis in equation (4.3) imply that all three groups on the right hand side of (4.5) are distinct. For ease of notation we write $H=W(\lk v)$. A free product such as (4.5) has a normal form (not to be confused with the normal form of a graph product) obtained as follows: choose complete sets of coset representatives for $W(v^*)/H$ and $W(\Gamma\setminus v)/H$. The former are the distinct elements of $\langle v\rangle$. As for $W(\Gamma\setminus v)/H$, we choose elements $b_{\nu}$ satisfying
 $$
 b_{\nu}\ \ \text{has no initial letter in}\ \ H.\tag4.6
 $$
Then [12], Theorem 4.4.1, each $g\in G$ has a unique normal form
$$
g=hv^{n_1}b_1\cdots v^{n_k}b_k,\tag4.7
$$
where $h\in H$ and only $h, n_1$ or $b_k$ may be trivial.
\proclaim{Lemma 4.5} The expression (4.7) is reduced in $G$.
\endproclaim
\demo{Proof} If $h\neq1$, $v\notin\supp h$ and $g=v^{n_1}hb_1v^{n_2}\cdots v^{n_k}b_k$ so that $v^{n_1}hb_1$ is reduced by (4.6). If $x$ is a last letter of $b_i$, $b_i\red b'x$ then $b_iv^{n_{i+1}}b_{i+1}=b'xv^{n_{i+1}}b_{i+1}$ is reduced, again by (4.6). If $x\in H$ $b_iv^{n_{i+1}}b_{i+1}=b'v^{n_{i+1}}xb_{i+1}$ is, again, reduced and $v$ is not a last letter of $b'$. If $x\notin H$, $b' v^{n_{i+1}}b_{i+1}$ is reduced as written.\qed
\enddemo
\proclaim{Lemma 4.6} Suppose $g$ is as in (4.7), $h'\in H$ and $h'g=gh'$. Then $h'$ commutes with each letter of $\supp b_i$, $1\leq i\leq k$ and with $h$. If $g$ has no initial letter in $H$, that is, $h=1$, then $h'$ commutes with each letter of $\supp g$.
\endproclaim
\demo{Proof} Since $gh'$ has normal form $h'g$, $b_kh'b^{-1}_k\in H$. Write $b_k=\beta_1\beta_2\beta_3\beta_4$ as in Corollary 4.2. Then $\beta_2, \beta_3\in H$ and $b_kh'b^{-1}_k=\beta_4\left(\beta_3\beta_2h'\beta^{-1}_2\beta^{-1}_3\right)\beta^{-1}_4=\beta_4h'' \beta^{-1}_4$ is reduced and it lies in $H$. Thus $\beta_4\in H$ contradicting (4.6). Thus $b_k=\beta_1$ and proceed by induction.\qed
\enddemo
By [13], Theorem 4.5, if $u\in G$ is CR and $uw=wu$, then one of the following occurs:
\roster
\item $ gug^{-1}\in H$, some $g\in G$,
\item$ gwg^{-1}\in H$, some $g\in G$,
\item both $gug^{-1}$ and $gwg^{-1}$ lie in the same factor of (5), some $g\in G$, or
\item there exist $g,z\in G$, $c,c'\in H$ so that $u=g^{-1}cgz^i$ and $w=x^{-1}c'gz^j$, for some $i,j$, and $g^{-1}cg$, $g^{-1}c'g$ and $z$ commute in pairs.
\endroster
Neither (1) nor (3) may occur in our case.
\proclaim{Lemma 4.7} If $u$ is CR, either $wu$ or $uw^{-1}$ is reduced and $\langle\supp u\rangle^c$ is connected, then $\langle u,w\rangle$ is cyclic.
\endproclaim
\demo{Proof} By Proposition 4.3, $\supp w=\supp u$. We analyze cases (2) and (4) above: if $w\in H$, that is if $g=1$, $\supp w\subset\lk v$ contradicting  (4.4). If $w\red g^{-1}h_0g$ we may assume that $g$ has no initial letter in $H$. If $wu$ is reduced, so is $uw$ and $g^{-1}ug$ is reduced and has no initial letter in
$H$. Say $h_0\red h^{-1}_1h'h_1$,  $h'$ CR. Then $h_1gug^{-1}h^{-1}_1$ is reduced with initial segment $h_1$. It follows that the normal form (4.7) has first factor $h_1$ and $h'$ commutes with $h_1$ and with the letters of $\supp b_i$. Thus $h_1$ must be trivial and $h_0$ is CR. By Lemma 4.6 $h_0$ commutes with each letter of $\supp u$. Write $g=g_4g_3g_2g_1$ as in Corollary 4.2. Then $gug^{-1}=g_4u'g^{-1}_4$ for some cyclic permutation of $u$ with the same support. Since $h_0$ commutes with $g_4u'g^{-1}_4$ and it is CR, Proposition 4.3 implies $\supp u\subset\supp{g_4u'g_4^{-1}}\subset H$ and this is a contradiction. This takes care of case (2). For (4), let $u=g^{-1}cgz^i$ and $w=g^{-1}c'gz^j$ and let $z_0=g^{-1}zg$ so that $u=g^{-1}cz_0g$ and $gug^{-1}\red g_4u'g^{-1}_4$ as above equals $cz^i_0$. It has no initial segment in $H$, since $g_4$ does not,  and so $z_0\red c^{-1}z_1$ and either $i=1$ or $c=1$.

If $c=1$ then $u=z^i$ so that $z$ is CR. Since $g^{-1}c'g$ commutes with $z$, the previous case implies $c'=1$ or $w=z^j$ and $\langle u,w\rangle\subset \langle z\rangle$.

If $i=1$, then $z=g^{-1}c^{-1}gu$ and $w=g^{-1}(c'c^{-j})gu^j$ and $g^{-1}(c'c^{-j})g$ commutes with $u$ so that $c'=c^j$, $w=c^jz^j=(cz)^j= u^j$ and again $\langle u,w\rangle$ is cyclic.\qed
\enddemo
The rest of the proof of the Centralizer Theorem follows as in [14].

\subhead The Centralizer Complex
\endsubhead We begin with a very important Lemma that we name the
\proclaim{Finite Order Lemma} Suppose $u\in G$ is CR, it has finite order and $\langle\supp u \rangle^c$ is connected. Then $\supp u$  consists of a single vertex. In particular, the order of $u$ lies in $P$. If $\theta\in\aut G$, $v\in\Gamma$, $\langle\supp\theta(v)\rangle^c$ is connected and $v$ has finite order then $\theta(v)$ is a conjugate of $v^m$,  $(m, o(v))=1$, so that $o(\theta(v))\in P$.
\endproclaim
\demo{Proof} All assertions follow from the first. If $o(u)<\infty$  and $u$ is CR, then $\supp u$ lies in a complete subgraph $\Delta\subset\Gamma$. Consequently $\langle\supp u\rangle^c\subset\Delta^c$ which is discrete. If it is also connected (in $\Gamma^c$) then $\supp u$ is a single vertex.\qed
\enddemo

 Suppose $\Gamma=(V,E)$ and $\Delta=(X,F)$ are two graphs. The join $\Gamma\ast \Delta$ is a graph with vertices $V\sqcup X$ and with edges, all the vertices of the individual graphs plus edges of the form $\{v,x\}$ for $v\in\Gamma,\  x\in\Delta$. Thus the individual graphs embed as subgraphs and any vertex of $\Gamma$ is adjacent to every vertex of $\Delta$. Then the centralizer of a CR element $u$ of $G=W(\Gamma,o)$ as expressed by (4.2) can be expressed as the graph product of a labeled graph $(K(u), o_K)$: if $\Lambda$ is a complete graph with vertices $ x_1\dots, x_n$ then $K=\Lambda\ast\lk\Sigma$. For the label we put $o_K(x_i)=o( u_i)$ and $o_K\restriction\lk\Sigma=o\restriction\lk\Sigma$. If $u_i$ has finite order then by the Finite Order Lemma, its support consists of only one point and its order lies in $P$. Now it is elementary that $W(K(u),o)$ is precisely formula (4.2). Thus by the Rigidity Lemma the element $u$ determines $(K(u),o)$, which has $n+|\lk\Sigma|$ vertices.
\definition{Definition 4.8} Suppose $g\in G$ and $g\red wuw^{-1}$ with $u$ CR. The rank $\rk g$ is the cardinality of the vertex set of $K(u)$.
\enddefinition

\head
5. $\aut G$
\endhead

As usual $\Gamma=(V,E)$ is a graph and $(\Gamma, o)$ is a labeled graph. We denote $W(\Gamma, o)$ by $G$. A mapping $\phi:V\to G$ is well-defined if it extends to an endomorphism of $G$, also denoted $\phi$. This means that if $\{x,y\}\in E$ then $\phi(x)$ and $\phi(y)$ commute and, additionally, $o(\phi(x))=o(x)$ for all $x\in\Gamma$.

\proclaim{Lemma 5.1} Labeled graph automorphisms, factor automorphisms and partial conjugations are always well-defined automorphisms.
\endproclaim

\definition{Definition 5.2} We say that an automorphism $\phi$ is conjugating if $\phi(x)$ is  a conjugate of $x$, for all $x\in\Gamma$. The subgroup of conjugating automorphisms of $G$ is denoted $\auto G$. More generally, given $\theta$, we write $\Cal C(\theta)$ to be the set of vertices $v$ for which $\theta(v)$ is a conjugate of $v$.
\enddefinition
Thus $\theta$ is conjugating precisely when $\Cal C(\theta)=\Gamma$; in general, $\Cal C(\theta)=\Cal C(\theta^{-1})$.
\proclaim{Proposition 5.3} The subgroup $\auto G$ is generated by the partial conjugations $\sigma_{K,v}$.
\endproclaim
A complete proof for the abelian case can be found in [10], Theorem 4.1. If all $o(x)=\infty$ (the right-angled Artin case) the result is also proved in [11], Theorem 2.2.
\definition{Definition 5.4} Suppose $x, y\in\Gamma$ we write
\roster
\item $x\leq y$ to mean $\lk x\subset y^{\star}$,
\item $x\leq_s y$ to mean $x^{\star}\subset y^{\star}$,
\item $x\sim y$ to mean $x\leq y$ and $y\leq x$,  and
\item $x\sim_s y$ to mean $x\leq_s y$ and $y\leq_s y$.
\endroster
Notice that $x\leq_s y$ if and only if $x$ and $y$ are adjacent and $x\leq y$. We write $[x]$ and $[x]_s$ for the equivalence classes with respect to $\sim$ and $\sim_s$, respectively.
\enddefinition
\proclaim{Proposition 5.5}  If $\tau=\tau_{x,y}$ is as above then it is well-defined if and only if
\roster
\item Either $o(x)=\infty$ and $x\leq y$ or,
\item $o(x)<\infty$, $x\leq_s y$ and $o(y)\ |\ o(x)$.
\endroster
If $x\leq_s y$,  $o(x)=p^j$, $o(y)=p^k$ and $k>j$, then $\tau$ is not well-defined but, if $q=p^{k-j}$,  $\tau^q$ is,  and $q$ is the smallest such positive power.
\endproclaim
\proclaim{Lemma 5.6} If $x\in\Gamma$, $[x]$ is either complete or discrete, $[x]_s$ is always complete and $[x]_s\subset [x]$. Equality holds if either $[x]_s$ has more than one element or if $[x]$ is complete.
\endproclaim
This is Proposition 3.2 of [11]. We now state Proposition 3.5 of [11] and to do so we recall Definition 4.8 and  introduce the important subgraph $\Gamma_v$ of $\Gamma$:
\definition{Definition 5.7} Let  $\Gamma_v:=\langle\{x\in\Gamma: v\leq x\}\rangle_\Gamma$. Similarly, let $\Omega_v=\Gamma_v\setminus\lk v$.
\enddefinition
\proclaim{Proposition 5.8} Suppose $g\in G$ is CR and $g=z_1^{r_1}\cdots z_s^{r_s}$ is a basic form with $Z_i^c=\langle\supp z_i\rangle^c$, connected. Let $Z=\bigsqcup Z_i=\supp g$; then, for all $y\in Z$, 
$$
\rk g\leq\rk y.
$$
Suppose that $z\in Z_1$ satisfies $\rk z =\rk g$. Then
\roster
\item For $i\geq2$, $z_i$ is a vertex,
\item $Z\subset \Gamma_z$;  in particular, $z\leq z_i$ for $i\geq2$, and
\item $Z_1\cap z^*=z$.
\endroster
\endproclaim
Thus, Proposition 5.8(3) asserts that $Z_1\subset\Omega_z$.
\proclaim{Lemma 5.9}If $\theta\in\aut G$, there exists an automorphism $\gamma\in\aut(\Gamma, o)$ that satisfies $v\in\csupp\gamma\theta(v)$ for all $v\in\Gamma$.
\endproclaim
This is an immediate consequence of the Rigidity Lemma. It also emphasizes the fact that certain properties of an automorphism $\theta\in\aut G$ will be used often. We remark on three of them:
$$
v\in\csupp\theta(v),\tag5.1
$$

$$
\langle\csupp\theta(v)\rangle^c\ \ \text{is connected.}\tag5.2
$$
and
$$
\csupp\theta(v)\subset \Gamma_v.\tag5.3
$$
\subhead Some Definitions
\endsubhead 
Many of the definitions below are temporary. First, let $E(G)$ (for "elementary") be the subgroup of $\aut G$ generated by the set $\Cal G$ of the Introduction. Our goal is to show that $E(G)=\aut G$. Second, we say that $\phi,\ \psi\in\aut G$ are E-equivalent if $\psi=\epsilon_1\phi\epsilon_2$ for some $\epsilon_i\in E(G)$. 
\definition{Definition 5.10} Let $\theta\in\aut G$; we say that $\theta$ is simple if it satisfies (5.1) and (5.2) for all $v\in\Gamma$. More generally, given $\theta$, we write $\Cal S(\theta)$ to be the set of 
vertices that satisfy (5.1) and (5.2).
\enddefinition
Thus, $\theta$ is simple precisely when $\Cal S(\theta)=\Gamma$.
Assume that $\chi$ is simple. Since $v\in\csupp\chi(v)$ and the rank of $v$ equals the rank of the CR part of $\chi(v)$, Proposition 5.8(2) shows that $\csupp\chi(v)\subset \Gamma_v$. This motivates the following
\definition{Definition 5.11} Let $\theta\in\aut G$; we say that $\theta$ is quasi-simple if it satisfies (5.2) and (5.3) for all $v\in\Gamma$.
\enddefinition
\proclaim{Lemma 5.12} Suppose $\theta\in\aut G$ and $\iota$ is an inner automorphism of G; then for all $v\in\Gamma$, $\csupp\theta(v)=\csupp\iota\theta(v)$. In particular, $\theta$ is simple (resp. quasi-simple) if and only if $\iota\theta$ is simple (resp. quasi-simple) and, for a chosen vertex $v$, we may assume that $\theta(v)$ is CR
\endproclaim 
\demo{Proof} This follows from the Reduction Lemma: if $\iota$ is conjugation by $w$ and $u$ is the CR part of $\theta(v)$, as in Corollary 4.2, $w\theta(v)w^{-1}=w_4 u' w_4^{-1}$ and $u'$ is a cyclic permutation of $u$ with the same support as $u$.\qed
\enddemo

\proclaim{Lemma 5.13} If $\chi$ is simple $\chi\langle\Omega_v\rangle=\langle\Omega_v\rangle$.
\endproclaim
This follows from Lemma 5.5 and its Corollary in [11].
\proclaim{Theorem 5.14}If $\chi$ is simple then $\chi^{-1}$ is quasi-simple.
\endproclaim
This is Theorem 5.6 of [11]. 

Recall that, if $\Sigma\subset\Gamma$, $\Sigma^*=\bigcap_{s\in\Sigma}s^*$  is the graph whose vertices  are adjacent to each vertex of $\Sigma$. 
The Lemma below is part of the proof of Theorem 6.3 in [11] and it requires a rather elaborate proof.
\proclaim{Lemma 5.15} Let $\theta\in\aut G$. Suppose that $\theta$  satisfies (5.1) for all vertices. If $\theta(v)$ is CR and $v\leq_s x$, then
$\theta(x)$ is also CR.
\endproclaim
\demo{Proof} We write $\theta(x)\red wuw^{-1}$ with $u $ CR and show that $w\neq1$ leads to a contradiction. We may assume $x \neq v$.  First apply Proposition 5.8 to $g=\theta(v)$: let
$\theta(v)=t_1^{r_1}\cdots t_s^{r_s}$ be in basic form; assume $T_i=\supp t_i$, $v\in T_1$ and $T=\bigcup_{i=1}^s T_i$. Since $\rk
v=\rk\theta(v)$ it follows that $t_2,\dots,t_s$ are vertices and that $T_1\cap v^\star=v$. We also have 
$$
(T\setminus T_1)\sqcup\lk T\subset \lk v.\tag 5.4
$$

Now apply Proposition 5.8 to $g=u$: let $u=u_1^{m_1}\cdots u_q^{m_q}$ be in basic form with $U_j=\supp u_j$, $x\in U_1$ and $U=\bigcup_j U_j$. Since $\rk x=\rk u$, the $u_j$, $j\geq2$ are vertices and
 by Proposition 5.8(2), for $y\in U$, we have $x\leq y$. Since $v\leq_s x$ by hypothesis and $x\leq y$, then $v\leq y$ and we have, in order,
 $$
v^\star\subset x^\star,\qquad v\in\lk x\subset y^\star\ \ \ \text{and}\ \ \ \lk v\subset y^\star     \tag5.5
 $$
 with $y\in U$. Consequently $v^\star\subset y^\star$ and 
 $$
 v^\star\subset U^\star.\tag5.6
 $$

Since $x\in v^\star$, $x$ commutes with every vertex of $U$, then, by Proposition 5.8(3), $U_1=x$. Moreover, since $v \in \lk x$, then (5.5) implies that  $v \in y^\star$. Therefore $y \in v^\star$, for all $y \in U$. So we have
$U \setminus v \subset \lk v$. This, together with the fact that  $T_1\cap \lk v= \emptyset$, imply that 

$\bullet$  if $v \notin U$, then $U \cap T_1 = \emptyset$.

Claim: if $v \in U$, then necessarily $U \cap T_1 = v$. Indeed, $v\in U=\supp u$ implies $x\leq v$ which in turn means that $v\sim_s x$. If $t\in T=\sup\theta(v)$, then $v\leq t$ and, by the equivalence, $x\leq t$.
Then $x\in\lk v$ and $\lk v\subset t^\star$ imply $v^\star=x^\star\subset t^\star$ and, as with (5.6)  $v^\star\subset T^\star$ so that $T_1=v$ proving the claim.

Now $x\in C(v)$ so that $\theta(x)\red wuw^{-1}\in C(\theta(v))$ and both $U=\supp u$ and $\sup w$ are in $T\sqcup \lk T$. Assume $w\neq1$ and let $l\in\supp w$. If $l=v$, then $l\in U^\star$ by (5.6); next assume $l\in 
T\setminus T_1$. Then (5.4) and (5.6) give $l\in U^\star$. Finally assume $l\in T_1\setminus v$. Then $T_1\neq v$ and, by the Claim, we necessarily have $U\cap T_1=\emptyset$. Then $U\subset (T\setminus T_1)\cup\lk T$ and 
$l$ commutes with each letter of the latter set so $l\in C(u)$. In conclusion, $w\in C(u)$ which contradicts the fact that $wuw^{-1}$ is reduced.\qed 
\enddemo

\proclaim{Corollary 5.16} If $\theta,\ v$ as above, $\theta(v)=t_1^{r_1}\cdots t_s^{r_s}$ in basic form with $v\in\supp t_1$. Then  $v\leq_s t_i$  for $i=2,\dots, s$ and $\theta(t_i)$ is CR.
\endproclaim
\proclaim{Lemma 5.17} If $\theta\in\aut G$ and $v$ is a vertex, then $\theta(\langle[v]_s\rangle)\subset\langle[v]_s\rangle$ implies that equality must hold.
\endproclaim
\demo{Proof} As in Proposition 6.1 of [11], recall that $\langle[v]_s\rangle$ is a finitely generated abelian group;  if $x\in[v]_s$ and $x\notin\theta(\langle[v]_s\rangle)$ then, if $o(x)=\infty$, an argument on the free part of the abelian group as in [11] gives a contradiction. If $o(x)=p^k$, a similar argument on the torsion subgroup generated by all vertices $y$ with $p\ |\ o(y)$ gives a contradiction. The rest follows as in [11].
\enddemo
\proclaim{Theorem 5.18} Every automorphism of $G$ is E-equivalent to a simple one.
\endproclaim
\demo{Proof} Our task is to take an automorphism $\theta$ and alter $\Cal S(\theta)$; now we show that for all vertices $v$, there exists a $\beta_v\in E(G)$ such that $v\in\Cal S(\theta\beta_v)$ and $\Cal S(\theta)\subset\Cal S(\theta\beta_v)$; this will complete the proof.
Suppose not and let $v$ be a counter-example of maximal rank. By Lemma 5.12 we may assume that $\theta(v)$ is CR and for every vertex $x$ with rank strictly larger than $\rk v$, we may assume that $x\in\Cal S(\theta)$ by maximality. Then we have that for all $\beta\in E(G)$ with $v\in\Cal S(\theta\beta)$ we must have $\Cal S(\theta)\nsubseteq\Cal S(\theta\beta)$. Write $\theta(v)$ in basic form as in Corollary 5.16;

Claim 1: in Corollary 5.16 we may assume that the $t_i\sim_s v$ for $i\geq2$.
Suppose not, we show that if, say, $t_2$ is not equivalent to $v$, we may eliminate it by composing with a transvection. Then, $\rk t_2>\rk v$ so $t_2\in\Cal S(\theta)$, by Corollary 5.16  $\theta(t_2)$ is CR and its support contains $t_2$. But $t_2$ also lies in $C(\theta(v))$ and in $v^*$, so $\theta(t_2)$ lies in $\theta(\langle v^*\rangle)\subset C(\theta(v))$ and by connectivity we get $\theta(t_2)\in\langle t_2\rangle$.
If $o(t_2)=\infty$ then $\theta(t_2)=t_2^{\pm1}$; say $o(t_2)=p^j$. Then $\theta(t_2)=t_2^m$ with $m$ mutually prime with $p$. If $o(v)=o(\theta(v))=q^k$, then $(t_1^{r_1}\cdots t_s^{r_s})^{q^k}=1$ and so $t_2^{r_2q^k}=1$ and we have
$$
p^j\ |\ r_2q^k,\tag5.7
$$
then $p=q$ for otherwise $p^j\ |\ r_2$ and $t_2^{r_2}=1$. If $j>k$, let $d=p^{j-k}$; then (5.7) implies $d\ |\ r_2$ so write $r_2=ds$ and choose $m'$ to be an inverse of $m$ modulo $p^j$. If $\tau=\tau_{v,t_2}$, then $\tau^d$ is well-defined and $\tau^{m'sd}(v)=vt_2^{-m'r_2}$. Applying $\theta$ to both sides we have
$$
\theta\tau^{m'sd}(v)=\theta(v)\theta(t_2^{-m'r_2})^m=\theta(v)t_2^{-r_2}=t_1^{r_1}t_3^{r_3}\cdots t_s^{r_s}.
$$
If $j\leq k$ or in the infinite case $\tau_{v,t_2}$ is defined and  we may compose $\theta$ with the transvection $\tau_{v,t_2}^{-r_2}$ and the composite maps $v$ to a product not containing $t_2$, as in the previous case. 

This proves the claim.

We know $s>1$ in Corollary 5.16, otherwise we would have (5.2) contrary to the hypothesis $v\notin\Cal S(\theta)$. Thus we have at least $t_2\sim_s v$. If $|T_1|\geq2$, let $y\neq v$ in $T_1$. By Proposition 5.8 $y$ and $v$ are not adjacent but $t_2y=yt_2$ contradicts $t_2\sim_s v$. Consequently, $T_1=v$ and either $t_1=v^{\pm1}$ or, when $o(v)<\infty$ equal to $v^m$, where $(m,o(v))=1$. Thus, $\langle v\rangle=\langle v^m\rangle$ and $\theta(v)\in\langle[v]_s\rangle$. The transvections used on Claim 1 fix all vertices except $v$. In that way, we may assume that $\theta(x)\in\langle[x]_s\rangle=\langle[v]_s\rangle$ for all $x\in [v]_s$ and, by Lemma 5.17 we have $\theta(\langle[v]_s\rangle)=\langle[v]_s\rangle$.

Claim 2. Define $\eta:V\to G$ by
$$
\eta(z)=\cases
\theta(z) & z\in [v]_s\\
z & z\notin[v]_s,
\endcases
$$
then $\eta$ is well-defined.

To see this notice that $[v]_s$ is complete. Take adjacent vertices $x\notin[v]_s$ and $y\in[v]_s$. Then $\theta(y)$ is a product of elements $t$ in $[v]_s$ and each $t^*=v^*$. Since $x\in y^*=v^*$, $x$ commutes with each $t$ and $\eta$ is well-defined. Then $\eta$ and its inverse are in $E(G)$ and $\theta\eta^{-1}\restriction[v]_s=\text{id.}$ so $v\in\Cal S(\theta\eta^{-1})\supset\Cal S(\theta)$, a contradiction.\qed
\enddemo

 \proclaim{Lemma 5.19} Suppose that $\theta$ satisfies (5.3) for all $v\in\Gamma$. Then $[v]\cap\csupp\theta(v)$ is not empty.
\endproclaim
\demo{Proof}By Lemma 5.9 there exists  a $\gamma\in\aut(\Gamma,o)$ satisfying $v\in\csupp\gamma\theta(v)=\gamma(\csupp\theta(v))$. Consequently $v=\gamma(x)$ for some $x\in\csupp\theta(v)$, and so $v\leq x$ by Equation (5.3). Say $v$ and $x$ are not adjacent; then $\lk v\subset\lk x$. On the other hand, $\gamma(\lk x)=\lk v$ and so they have the same (finite) cardinality and equality must hold. If $v, x$ are adjacent, a similar argument with the stars shows $v^\star=x^\star$. In either case $v\sim x$.\qed
\enddemo
This is part (3) of the definition of quasi-simple in [11], Definition 5.1.
\proclaim{Proposition 5.20} Suppose $x,y$ and $z$ are three vertices $x,y$ in the same component of $\Gamma\setminus z^\star$. If $\chi$ is quasi-simple, then $z\in\supp\chi(x)$ if and only if  $z\in\chi(y)$.
\endproclaim
\demo{Proof} By Lemma 5.12 we may assume that $\chi(x)$ is CR. We may also assume that $x$ and $y$ are adjacent. Let then $\chi(x)=u$ and $\chi(y)\red wu' w^{-1}$ with both $u$ and $u'$ CR.
Then $u\in wC(u')w^{-1}$. Now $u'$ has a connected support in $\Gamma^c$ and it has no proper roots so the Centralizer Theorem implies that $C(u')=\langle u'\rangle\times \langle\lk u'\rangle$ and so $u$ lies in                            
$w\langle u'\rangle w^{-1}\times w\langle\lk u'\rangle w^{-1}$. Consequently, we may write $u=w(u')^n w^{-1}\cdot wcw^{-1}=\chi(y)^n\cdot wcw^{-1}$, where $n$ is an integer and $c\in\langle\lk             u'\rangle$. Say $z\in\supp u$; then either 
$z\in\supp\chi(y)$ or $z\in\supp (wcw^{-1})$. If $z\in\supp w$ then $z\in\supp\chi(y)$. If $z\in\supp c$ then $z$ commutes with each letter of $\supp u'$. By Lemma 5.19 there is a vertex $v\in\supp u'$ with  $v\sim y$. If $zv=vz$ then necessarily $zy=yz$ which contradicts the fact that $y\in\Gamma\setminus z^\star$.\qed
\enddemo
\proclaim{Theorem 5.21} Every simple automorphism $\chi$ is E-equivalent to a conjugating automorphism in $\auto G$.
\endproclaim
\demo{Proof} Recall Definition 5.2: $\Cal C(\chi)=\{x\in\Gamma: \chi(x)=wxw^{-1}\}$. The strategy here is the same as in Theorem 5.18, namely we show if $v\notin\Cal C(\chi)$ then we have $\alpha\in E(G)$ with $v\in\Cal C(\chi\alpha)\supset\Cal C(\chi)$. The result then follows by induction on $|V|$. First, we may assume that $v\in\Cal C(\chi)$ if $o(v)<\infty$ by the Finite Order Lemma. Thus we may assume from the start that $v$ is of infinite order and of maximal rank among such vertices. Consequently, if $x\in \Gamma_v$ is not equivalent to $v$, we may assume that $x\in\Cal C(\chi)$ by maximality. Finally, by Lemma 5.12, we may assume that $\chi(v)$ is CR.

By Proposition 5.8(3), $\supp\chi(v)\subset\Omega_v$. Assume that the connected components of $\Gamma\setminus v^\star$ are $K_1,\dots,K_m$. Write $\Omega_i=K_i\cap\  \Omega_v$ and let these intersections be non-empty only for $1\leq i\leq n$. If $x\in\Omega_v\cap [v]$ and $x\in\Omega_i$, then $\Omega_i=x$ because any path out of $x$ must go through $\lk x=\lk v$. So we write $\Omega_v\cap [v]=\{v'_1,\dots, v'_r\}\cup\{v_1,\dots, v_t\}$, where the $v'_j$ are of finite order and the $v_k$ are of infinite order, with $v=v_1$. In conclusion, we may write 
$$
\Omega_v=\Omega_1\sqcup\cdots\sqcup\Omega_s\sqcup \{v'_1,\dots, v'_r\}\sqcup\{v_1,\dots, v_t\},
$$
and $s+r+t=n$.

Let $H=\langle\bigcup\Omega_i\cup\{v'_j\}\rangle$ and $F=\langle\{v_k\}\rangle$. The latter is a free group and $H$ is the free product of the $\langle\Omega_i\rangle$ and the finite cyclic groups $\langle v'_j\rangle$. Then $\langle\Omega_v\rangle=H\ast F$ and by Lemma 5.15 $\chi(H\ast F)=H\ast F$.

Claim 1. If  $z\in\Omega_v\setminus\Omega_i$, define $\sigma: V\to G$ by $\sigma(x)=zxz^{-1}$ if $x\in\Omega_i$ and the identity on $V\setminus\Omega_i$. Then $\sigma$ is well-defined. Assume  $x\in\Omega_i$ and $y\in\Omega_v\setminus\Omega_i$ and neither vertex commutes with a third vertex $z$, then $x$ and $y$ lie in distinct components of $\Gamma\setminus z^\star$ by  Proposition 6.5(2) in [11].  Then  $\sigma$ is well defined because it extends to $\sigma_{K,z}$, where $K$ is the component of $\Gamma\setminus z^\star$ containing $x$.

Claim 2. For $i\leq s$ $\chi(\langle\Omega_i\rangle)=w\langle\Omega_i\rangle w^{-1}$, for some $w\in G$: say $x\in\Omega_i$, $\chi(x)=wxw^{-1}$ and $\iota$ is the inner automorphism given by conjugation by $w^{-1}$. Then $\iota\chi(x)=x$ and $\iota\chi$ is simple. If $y\in\Omega_i$ and $z\in\Omega_v\setminus\Omega_i$, then $x,y$ lie in the same component of $\Gamma\setminus z^\star$ by Proposition 6.5(1) of [11] and by Proposition 5.21 we know that $z\notin\supp\iota\chi(y)$ so that $\iota\chi(\langle\Omega_i\rangle)\subset\langle\Omega_i\rangle$. Since $\chi^{-1}\iota^{-1}$ is quasi-simple by Theorem 5.16, equality holds and $\chi(\langle\Omega_i\rangle)=w\langle\Omega_i\rangle w^{-1}$.

Following [4], page 116, the generators of $\aut (H\ast F)$ are (i) allowable permutations, (ii) factor automorphisms and (iii) Whitehead automorphisms. These are generators by [3]. Equation (5.1) allows us to assume that $\chi\restriction H\ast F$ has no factors of type (i). We say that $\omega\in\aut (H\ast F)$ is a Whitehead automorphism if there is a non-trivial element $a$ in one of the factors of $H$ or of the form $v_k^{\pm1}$ so that $\omega$ either conjugates a factor of $H$ by $a$ or leaves it fixed, pointwise in all cases, and $\omega(v_k)$ is one of $v_k$, $av_k$, $v_ka^{-1}$ or $av_ka^{-1}$. If one of the  $\langle\Omega_i\rangle$ contains an infinite cyclic free factor $\langle x\rangle$,  then we also have Whitehead automorphisms $\omega$ that map $x$ to one of $x$, $ax$, $xa^{-1}$ or $axa^{-1}$; however, Claim 2 implies that in that case only conjugations by $a$ are allowed. With the help of factor automorphisms we can consider three possibilities for $\omega(v_k)$: $v_k$, $v_ka$ and $av_ka^{-1}$. 

Relation 6 on [4], page 117, allows us to write $\chi\restriction H\ast F=\beta\alpha'$, where $\beta$ is a product of factor automorphisms fixing $F$. Then $\beta$ extends to an automorphism (also denoted $\beta$) in $E(G)$. The remaining factors, in $\alpha'$, are the inversions in $F$, transvections $\tau_{v_k,a}$ and conjugations of one factor of $H$ by $a$. The first two types extend to all of $G$ by mapping all vertices in $V\setminus v_k$ to themselves and the conjugations extend by Claim 1. Let $\alpha'$ the automorphism of $E(G)$ so obtained and denote by $\alpha$ its inverse.

Then $\chi\alpha$ maps each $v_k$ to itself and $\Cal C(\chi)\subset\Cal C(\chi\alpha)$, and we are done.\qed

\enddemo
The Main Theorem follows from Theorem 5.22 and Proposition 5.3.

\head
6. $\aut^\star G$
\endhead
In this section we concentrate on an important subgroup of $\aut G$ defined in [5]. To that effect we let $MCS(\Gamma)$ be the set of all maximal complete subgraphs (or cliques) of $\Gamma$. Then we define
\definition{Definition 6.1} The subgroup $\aut^\star G$ of $\aut G$ is the subgroup of all $\theta$ such that, for all $\Delta\in MCS(\Gamma)$,  $\theta(W(\Delta))$ is a conjugate of $W(\Theta)$ for some $\Theta$ in $MCS(\Gamma)$. The subgroup $\aut^1 G$ is made up of all $\theta$ with $\theta(W(\Delta))= W(\Theta)$.
\enddefinition
By Theorem 3.1 of [7], $\aut^\star G$ is the semi-direct product $\auto G\rtimes \aut^1 G$. Here we intend to find a set of generators for $\aut^1 G$.
\definition{Definition 6.2} Let $\tau_{x,y}$ or $\tau_{x,y}^q$ be a well-defined transvection in $\Cal G$. We say that the transvection is a $\star$-transvection if $x$ is adjacent to $y$. The set $\Cal G^\star$ is the set of all labeled graph automorphisms, factor automorphisms, $\star$-transvections and partial conjugations. The set $\Cal G^1$ is obtained from $\Cal G^\star$ by removing all partial conjugations.\enddefinition
Our aim is to show
\proclaim{Theorem 6.3} $\aut^\star G$ is the subgroup of $\aut G$ generated by $\Cal G^\star$.
\endproclaim
\demo{Proof }  By the results of [7], it suffices to show that $\Cal G^1$ generates $\aut^1 G$. Suppose $\eta\in\aut^1 G$; then, for some $\gamma\in\aut(\Gamma,o)$ we have that $\theta=\gamma\eta\in\aut^1 G$ satisfies (5.1) for all $v\in\Gamma$. Let $\Delta, \Theta\in MCS(\Gamma)$ as in Definition 6.1 and suppose $\Delta$ has vertices $v$. Then all elements of $\theta(\langle\Delta\rangle) =\langle\Theta\rangle$ are CR (see the remarks in [7], \S 3.1) and if $u=\theta(v)$,  we may write it in basic form as $u=u_1^{r_1} \cdots u_s^{r_s}$ and $C(\theta(v))\subset C(u_i)$. By the maximality of $\langle\theta(\Delta)\rangle$ the $u_i\in\theta(\langle\Delta\rangle)$.  Consequently, $\theta(\langle\Delta\rangle)$ is a product of cyclic groups of the form $\langle u_i\rangle$, with $\langle\supp u_i\rangle^c$ connected (See Proposition 8.3 and Lemma 8.5 of [10]). Each $v\in\Delta$ must lie in a $\supp u_i$ by (5.1) and if $\theta(\Delta)\cap\supp u_i=\emptyset$, then $\langle\theta(\Delta),u_i\rangle$ is an abelian group which contains $\langle\theta(\Delta)\rangle$ and that contradicts maximality. Since $\Theta^c$ is discrete, it follows that each $\supp u_i$ consists of a single vertex of $\Delta$ and by maximality $\theta(\langle\Delta\rangle)=\langle\Delta\rangle$.
If $\Lambda\in MCS(\Gamma)$ then $\theta$ must map $\langle\Delta\cap\Lambda\rangle$ to itself (see [5]). Then $\theta_{\Delta}:=\theta\restriction \langle\Delta\rangle$ extends to $G$ by the formula $\theta_{\Delta}(z)=z$ if $z\notin\Delta$. The results of [9] show that $\theta_{\Delta}$ is a product of $\star$-transvections and factor automorphisms (the elementary automorphisms of an abelian group) and $\theta
\theta_{\Delta}^{-1}\in\aut^1 G$ fixes $\langle\Delta\rangle$ and we proceed by induction on $|MCS(\Gamma)|$.\qed
\enddemo
\remark{Remarks 6.4} For simplicity we assume that $G$ is a right-angled Artin group.
\roster
\item Let $N$ be the subgroup of $\aut^1 G$ generated by the factor automorphisms $\phi_x$ and the $\star$-transvections. If $\gamma\in\aut (\Gamma,o)$ 
$\gamma\phi_x\gamma^{-1}=\phi_{\gamma(x)}$ and a similar formula holds for transvections. Therefore $N$ is normal. Unfortunately, in general $N\cap\aut(\Gamma,o)\neq 1$ as the example below shows.
\item It is possible to express elements of $\aut (\Gamma,o)$ as products of $\star$-transvections and factor automorphisms. If $\Gamma$ has four vertices $a,b,c,d$ and its edges are $\{\{a,b\}, \{a,c\}$,$\{b,c\}, \{c,d\}\}$ the only non-trivial graph automorphism is the flip of $a$ and $b$. But this can be accomplished with $\tau_{a,b},\tau_{b,a}$ and some inversions $\phi_x$.
\item Similarly, it is possible to express elements of $\auto G$ as products of non-$\star$-transvections and inversions: if $x\leq y$ but $x,y$ are not adjacent, then $x$ is a connected component of $\Gamma\setminus y^\star$ so that $\sigma_{x,y}$ is defined. But an easy calculation gives  $\sigma_{x,y}=(\tau_{x,y}^{-1}\phi_x)^2$. There are also many instances in which $\aut G=\aut^\star G$, for example if the graph is a hexagon (cf. [7], Lemma 1.2).
\endroster
\endremark
\head 
7. Whitehead Automorphisms
\endhead
It is interesting to speculate whether the work in [2] can be generalized to arbitrary graph products of finitely generated abelian groups. We conjecture that the answer is yes. A first step would be to find a good set of generators $\aut G$, and experience shows that $\Cal G$ is not the best set. Here we quickly give a sketch of a description of the Whitehead automorphisms as defined in [2], \S3 and [4],  \S2. Labeled graph automorphisms and factor automorphisms are termed Whitehead automorphisms of type I and we call the set of these $\Omega_1$ Observe that $\langle\Omega_1\rangle$ is a finite subgroup of $\aut G$. To define Whitehead automorphisms of type II, which form a second set $\Omega_2$, we let $S$ be the set of all vertices  $v\in\Gamma$ with $o(v)=\infty$ and   write $J$ for the set of vertices $v\in\Gamma$ with $o(v)<\infty$. Then let $L$ be $S\sqcup S^{-1}\sqcup\bigsqcup_{v\in J} (\langle v\rangle\setminus1)$. Let $a$ in $L$. Then $a=v^q$ for some $q$ if $v\in J$ or $s^{\pm1}$ for $s\in S$.
\definition{Definition 7.1} Given $a$ as above, define $\phi$ as follows:
\roster
\item if $a=v^q$ and $x\in J\setminus v^*$, then $\phi(x)$ is either $x$ or $axa^{-1}$, 
\item if $a=v^q$ and  $x\in v^*$, then $\phi(x)$ is one of $x$, $axa^{-1}$, $a^mx$, $xa^{-m}$, with $m$ minimal for $o(a^m)\ |\ o(x)$
\item if $a=s^{\pm1}$ and $x\in J$, then $\phi(x)$ is either $x$ or $axa^{-1}$, and
\item if $a=s^{\pm1}$ and $x\in S\sqcup S^{-1}$, then $\phi(x)$ is one of $x$, $axa^{-1}$, $ax$ or $xa^{-1}$.
\endroster
\enddefinition 
Then we define a set $A\subset L$ by
$$
\{v\in J: \phi(v)=ava^{-1}\  \text{or} \  \ a^mv\}\cup\{s\in S\sqcup S^{-1}: \phi(s)= asa^{-1}\ \text{or}\ \ as\}\cup\{\zeta\},
$$
where $\zeta=v$ if $a=v^q$ or $s$ if $a=s^{\pm1}$. See [3], \S2.
Then $(A,a)$ determine $\phi$. This is not always well-defined, but if it is, we refer to it as a type II Whitehead automorphism denoted by $(A,a)$. One can get the typical table of values for $(A,a)x$. For example, $(A,a)x=xa^{-m}$ if $x\notin A$, $x^{-1}\in A$ and $o(x)<\infty$; for $o(a)=\infty$ the behavior of $(A,a)$ is identical to that described on page 31 of [12].
\proclaim{Lemma 7.2} Suppose $A\subset L$ satisfies that if $a\in S\sqcup S^{-1}$ then $a\in A$ but $a^{-1}\notin A$. Then $(A,a)$ is well defined if the following hold:
\roster
\item  The set $V\cap A\cap A^{-1}\setminus \lk\zeta$ is a union of connected components of $\Gamma\setminus\zeta^\star$,
\item For each $x\in A\setminus A^{-1}$, $\lk x\subset\zeta^\star$, and
\item If $o(x)<\infty$ in (2), $x^\star\subset\zeta^\star$ and $o(\zeta^m)\ |\ o(x)$.
\endroster
\endproclaim
This is a generalization of  Lemma 2.5 in [2]. We now let $\Omega=\Omega_1\cup\Omega_2$. Then $\Cal G\subset\Omega$ and so the Whitehead automorphisms generate $\aut G$. 
An obvious modification on Definition 7.1(4) gives a set of generators for $\aut^\star G$.

\enddocument